\documentclass[envcountsame,oribibl]{llncs}
\usepackage{amsmath,amsfonts}

\newcommand{\rk}{{\operatorname{rk}}}
\newcommand{\den}{{\operatorname{den}}}
\newcommand{\num}{{\operatorname{num}}}
\newcommand{\smooth}{{\operatorname{smooth}}}
\newcommand{\FF}{{\mathcal F}}

\newcommand{\C}{{\mathbf C}}
\newcommand{\Q}{{\mathbf Q}}

\newcommand{\Z}{{\mathbf Z}}

\newcommand{\R}{{\mathbf R}}

\newcommand{\OO}{{\mathcal O}}
\newcommand{\pp}{{\mathfrak p}}

\newcommand{\F}{{\mathbf F}}

\newcommand{\Proj}{\operatorname{Proj}}

\newcommand{\nichts}{\ensuremath{\left.\right.}}
\newcommand{\isom}{\simeq}

\newcommand{\tors}{{\operatorname{tors}}}
\newcommand{\red}{{\operatorname{red}}}

\begin{document}

\title{Using elliptic curves of rank one towards the undecidability of Hilbert's Tenth Problem over rings of algebraic integers}
\titlerunning{Hilbert's Tenth Problem over rings of algebraic integers}
\author{Bjorn Poonen\thanks{This research was supported by 
NSF grant DMS-9801104, and a Packard Fellowship.
This article will appear in the Proceedings of the Algorithmic Number Theory
Symposium~V, in the Springer Lecture Notes in Computer Science series,
accessible electronically at \mbox{{\tt http://www.springer.de/comp/lncs/index.html}}  \copyright Springer-Verlag.}}
\institute{Department of Mathematics, University of California, Berkeley, CA 94720-3840, USA, \email{poonen@math.berkeley.edu}}

\maketitle

\begin{abstract}
Let $F \subseteq K$ be number fields,
and let $\OO_F$ and $\OO_K$ be their rings of integers.
If there exists an elliptic curve $E$ over $F$ 
such that $\rk\, E(F) = \rk\, E(K) = 1$,
then there exists a diophantine definition of $\OO_F$ over $\OO_K$.
\end{abstract}

\section{Introduction}
\label{introduction}

D.~Hilbert asked, as Problem~10 of his famous list of 23 problems 
posed to the mathematical community in 1900,
for an algorithm to decide, given a polynomial equation
$f(x_1,\dots,x_n)=0$ with coefficients in the ring $\Z$ of integers,
whether there exists a solution with $x_1,\dots,x_n \in \Z$.
In Hilbert's time, there was no formal definition of algorithm,
but presumably what he had in mind was a mechanical procedure
that a human could in principle carry out, 
given sufficient paper, pencils, erasers, and time, 
following a set of strict rules requiring 
no insight or ingenuity on the part of the human.
In the 1930s, several rigorous models of computation were proposed
as a substitute for the informal notion of ``mechanical procedure'' as above
(the $\lambda$-definable functions of A.~Church and S.~Kleene, 
the recursive functions of K.~G\"odel and J.~Herbrand, and the 
logical computing machines of A.~Turing).
These models, as well as others developed later, were shown to be equivalent; 
this gave credence to the {\em Church-Turing thesis},
which is the belief that every mechanical procedure 
can be carried out by a Turing machine.
Therefore, the modern interpretation of Hilbert's Tenth Problem
is that it asks whether a Turing machine can decide the existence of
solutions.

J.~Matijasevi{\v{c}}~\cite{matijasevic}, 
building on earlier work by 
M.~Davis, H.~Putnam, and J.~Robinson~\cite{davis-putnam-robinson1961}
showed that there is no such Turing machine.
To describe their work in more detail,
we need a few definitions.
A subset $S$ of $\Z^n$ is called {\em listable} or {\em recursively enumerable}
if there is an algorithm (Turing machine) such that $S$ is exactly
the set of $a \in \Z^n$ that are eventually printed by the algorithm.
A subset $S$ of $\Z^n$ is said to be {\em diophantine}, or to
admit a {\em diophantine definition},
if there is a polynomial 
$p(a_1,\dots,a_n,x_1,\dots,x_m) \in \Z[a_1,\dots,a_n,x_1,\dots,x_m]$
such that 
$$S = \{\, a \in \Z^n : (\exists x_1,\dots,x_m \in \Z ) \; 
	p(a_1,\dots,a_n,x_1,\dots,x_m)=0\,\}.$$
For example, the subset $\Z_{\ge 0}:=\{0,1,2,\dots\}$ 
of $\Z$ is diophantine,
since for $a \in \Z$, we have
	$$a \in \Z_{\ge 0} \iff (\exists x_1,x_2,x_3,x_4 \in \Z) \; 
	x_1^2+x_2^2+x_3^2+x_4^2 = a.$$

One can show using ``diagonal arguments''
that there exists a listable subset $L$ of $\Z$ whose complement
is not listable.
It follows that for this $L$, there is no algorithm that takes as input 
an integer $a$ and decides in a finite amount of time 
whether $a$ belongs to $L$;
in other words, membership in $L$ is undecidable.

Diophantine subsets of $\Z^n$ are listable: given $p$,
one can write a computer program with an outer loop
with $B$ running through $1,2,\dots$,
and an inner loop in which one tests the finitely many
$(a_1,\dots,a_n,x_1,\dots,x_m) \in \Z^{n+m}$
satisfying $|a_i|,|x_j| \le B$ for all $i$ and $j$,
and prints $(a_1,\dots,a_n)$ if $p(a_1,\dots,a_n,x_1,\dots,x_m)=0$.
Davis~\cite{davis1953} 
conjectured conversely that all listable subsets of $\Z^n$
were diophantine, and this is what Matijasevi\v{c} 
eventually proved.
In particular, the set $L$ is diophantine.
Hence a positive answer to Hilbert's Tenth Problem
would imply that membership in $L$ is decidable.
But membership in $L$ is undecidable,
so Hilbert's Tenth Problem is undecidable too;
that is, there is no algorithm that takes as input a polynomial
$p \in \Z[x_1,\dots,x_n]$,
and decides whether $p(x_1,\dots,x_n)=0$ has a solution in integers.

More generally, if $R$ is any commutative ring with $1$,
one can define what it means for
a subset of $R^n$ to be {\em diophantine over $R$},
by replacing $\Z$ by $R$ everywhere.
Similarly one can speak of {\em Hilbert's Tenth Problem over $R$}
provided that one has fixed some encoding of elements of $R$
as finite strings of symbols from a finite alphabet, 
so that polynomials over $R$
can be the input to a Turing machine.
For some rings $R$ (for example, uncountable rings)
such an encoding may not be possible.
In this case one should modify the problem,
by specifying a countable subset $\mathcal P$ of the set of all polynomials
over $R$ and an encoding of elements of $\mathcal P$ as finite strings
of symbols, and then asking whether there exists a Turing machine that
takes as input a polynomial $p \in \mathcal P$
and decides whether $p(x_1,\dots,x_n)=0$ has a solution over $R$.
For example, K.~Kim and F.~Roush~\cite{kim-roush1992}
proved that Hilbert's Tenth Problem over the purely transcendental
function field $\C(t_1,t_2)$ is undecidable
when one takes $\mathcal P$ to be the set of polynomials
with coefficients in $\Z[t_1,t_2]$.
Usually it is not necessary 
to specify exactly how the elements of $\mathcal P$
are encoded, since usually given any two reasonable encodings,
a Turing machine can convert between the two.

Perhaps the most important unsolved question in this area
is Hilbert's Tenth Problem over the field $\Q$ of rational numbers.
The majority view seems to be that it should be undecidable.
To prove this, it would suffice to show that the
subset $\Z$ of $\Q$ is diophantine over $\Q$.
On the other hand, B.~Mazur has suggested that perhaps
for any variety $X$ over $\Q$, 
the topological closure of $X(\Q)$ in $X(\R)$ has
at most finitely many connected components;
if this is true, no such diophantine definition of $\Z$ over $\Q$ exists.
See~\cite{mazur1994} and the 
more recent articles~\cite{cornelissen-zahidi2000}
and~\cite{pheidas2000} for further discussion.

The function field analogue, namely Hilbert's Tenth Problem over
the function field $k$ of a curve over a finite field,
is known to be undecidable.
The first result of this type
is due to T.~Pheidas~\cite{pheidas1991},
who proved this for $k=\F_q(t)$ with $q$ odd.
His argument was adapted and generalized by
C.~Videla~\cite{videla1994} for $k=\F_q(t)$ with $q$ even,
by A.~Shlapentokh~\cite{shlapentokh1992functionfield}
for other function fields of odd characteristic,
and finally by K.~Eisentr\"ager~\cite{eisentraeger2003}
for the remaining function fields of characteristic~2.
Analogues are known also for many function fields over infinite
fields of positive characteristic: see~\cite{shlapentokh2000functionfield}
and~\cite{eisentraeger2003}.

For more results concerning Hilbert's Tenth Problem,
see the book~\cite{hilbertstenthproblem},
and especially the survey articles~\cite{pheidas-zahidi2000}
and~\cite{shlapentokh2000survey} therein.
Since the publication of that book, 
undecidability of Hilbert's Tenth Problem has been proved
also for function fields of curves $X$ over formally real fields $k_0$
with $X(k_0)$ nonempty~\cite{moretbailly2002} 
(in fact this is just one application of his results),
and for function fields of surfaces over real closed or
algebraically closed fields
of characteristic zero~\cite{eisentraeger2003}.

\section{Hilbert's Tenth Problem over rings of integers}
\label{h10overOk}

In this article, our goal is to prove
a result towards Hilbert's Tenth Problem over rings of integers.
If $F$ is a number field, let $\OO_F$ denote the integral
closure of $\Z$ in $F$.
There is a known diophantine definition of $\Z$ over $\OO_F$
for the following number fields:
\begin{enumerate}
\item $F$ is totally real~\cite{denef1980}.
\item $F$ is a quadratic extension of a 
totally real number field~\cite{denef-lipshitz1978}.
\item $F$ has exactly one conjugate pair of 
nonreal embeddings~\cite{pheidas1988},~\cite{shlapentokh1989}.
\end{enumerate}
In particular, Hilbert's Tenth Problem over $\OO_F$ is undecidable
for such fields $F$.

It is conjectured~\cite{denef-lipshitz1978} 
that for {\em every} number field $F$,
there is a diophantine definition of $\Z$ over $\OO_F$.
Our main theorem gives evidence for this conjecture,
by reducing it to a plausible conjecture 
about the existence of certain elliptic curves.

Before stating our theorem,
let us recall the Mordell-Weil Theorem, which states
that if $E$ is an elliptic curve over a number field $F$, 
then the abelian group $E(F)$ is finitely generated.
Let $\rk\, E(F)$ denote the rank of $E(F)$.

\begin{theorem}
\label{main}
Let $F \subseteq K$ be number fields,
and let $\OO_F$ and $\OO_K$ be their rings of integers.
Suppose that there exists an elliptic curve $E$ over $F$ 
such that $\rk\, E(F) = \rk\, E(K) = 1$.
Then there exists a diophantine definition of $\OO_F$ over $\OO_K$.
\end{theorem}

Most of the rest of this paper
is devoted to the proof of Theorem~\ref{main}.
But for now, we mention its application to Hilbert's Tenth Problem.

\begin{corollary}
Under the hypotheses of Theorem~\ref{main},
if in addition $F$ is of one of the types of number fields listed
above for which a diophantine definition of $\Z$ over $\OO_F$
is known, then Hilbert's Tenth Problem over $\OO_K$ is undecidable.
\end{corollary}

\begin{proof}
Theorem~\ref{main} reduces the undecidability over $\OO_K$
to the undecidability over $\OO_F$.
\qed
\end{proof}

J.~Denef, at the end of~\cite{denef1980}, sketches a simple proof
of Theorem~\ref{main} in the case where $K$ is totally real
and $F=\Q$.
In fact, he is also able to treat some 
totally real algebraic extensions $K$
of {\em infinite} degree over $\Q$.
But his proof technique does not seem to generalize easily
to fields that are not totally real.

Our proof of Theorem~\ref{main}
is similar to that of an older result,
the theorem of~\cite{denef-lipshitz1978},
which uses a $1$-dimensional torus (a Pell equation)
in place of the elliptic curve.
We have been inspired also by the exposition 
of the ``weak version of the vertical method'' 
in~\cite{shlapentokh2000survey}
and by the ideas in~\cite{pheidas2000}.

\subsection{Preliminaries on diophantine sets over $\OO_K$}

The subset $\OO_K - \{0\}$ of $\OO_K$ is diophantine
over $\OO_K$: see Proposition~1(c) of~\cite{denef-lipshitz1978}.
We have a surjective map 
$\OO_K \times (\OO_K - \{0\}) \rightarrow K$ taking $(a,b)$ to $a/b$.
If $S \subseteq K^{n}$ is diophantine over $K$,
then the inverse image of $S$ under 
$(\OO_K \times (\OO_K - \{0\}))^{n} \rightarrow K^{n}$
is diophantine over $\OO_K$.
In this case, we will also say that $S$ is diophantine over $\OO_K$.
It follows that in constructing diophantine definitions over $\OO_K$,
there is no harm in using equations with some variables taking
values in $\OO_K$ and other variables taking values in $K$.

Given $t \in K^\times$, define the {\em denominator ideal} 
$\den(t) = \{\, b \in \OO_K :  bt \in \OO_K \,\}$
and the {\em numerator ideal} $\num(t) = \den(t^{-1})$.
Also define $\num(0)$ to be the zero ideal.
These ideals behave in the obvious way upon extension of the field.

\begin{lemma}
\label{ideals}
\nichts
\begin{enumerate}
\item For fixed $m,n \in \Z_{\ge 0}$,
the set of $(x_1,\dots,x_m,y_1,\dots,y_n)$ in $K^{m+n}$
such that the fractional ideal $(x_1,\dots,x_m)$
divides the fractional ideal $(y_1,\dots,y_n)$
is diophantine over $\OO_K$.
\item The set of $(t,u) \in K^\times \times K^\times$ such that
$\den(t)  \mid  \den(u)$ 
is diophantine over $\OO_K$.
\item The set of $(t,u) \in K^\times \times K$ such that
$\den(t)  \mid  \num(u)$ 
is diophantine over $\OO_K$.
\item The set of $(t,u) \in \OO_K \times K^\times$ such that
$t  \mid  \den(u)$ 
is diophantine over $\OO_K$.
\end{enumerate}
\end{lemma}

\begin{proof}
Statement~1 is clear, since the condition is that 
there exist $c_{ij} \in \OO_K$ such that $y_j=\sum_i c_{ij} x_i$ 
for each $j$.
Statement~2 follows from statement~1, since
\mbox{$\den(t) \mid \den(u)$} if and only if the fractional ideal $(u,1)$
divides $(t,1)$.
Statements 3 and~4 follow from statement~2: namely,
\begin{align*}
	\den(t) \mid \num(u) &\iff u=0 \text{ or } 
	(\exists v)(uv=1 \text{ and } \den(t) \mid \den(v)),\\
	t \mid \den(u) &\iff 
	(\exists v)(tv=1 \text{ and } \den(v) \mid \den(u)).
\end{align*}
\qed
\end{proof}

\subsection{Bounds from divisibility in $\OO_K$}

Let $n=[K:\Q]$ and $s=[K:F]$.
Fix $\alpha \in \OO_K$ such that $\{1,\alpha,\dots,\alpha^{s-1}\}$
is a basis for $K$ over $F$.
Let $D \in \OO_F$ denote the discriminant of this basis.
If $I$ is an ideal in $\OO_K$,
let $N_{K/\Q}(I) \in \Z_{\ge 0}$ denote its norm.

\begin{lemma}
\label{divisibilitybound}
There is a positive integer $c>0$ depending only on $F$, $K$, and $\alpha$
such that the following holds.
Let $I \subset \OO_K$ be a nonzero nonunit ideal,
and let $\mu \in \OO_K$.
Write $\mu = \sum_{i=0}^{s-1} a_i \alpha^i$ with $a_i \in F$.
Suppose that $\mu(\mu+1)\cdots(\mu+n) \mid I$.
Then $N_{K/\Q}(Da_i) \le N_{K/\Q}(I)^c$.
\end{lemma}

\begin{proof}
This is essentially Section~1.2 of~\cite{shlapentokh2000survey}.
The only differences are that we have specialized by taking $l_i=-i$,
and we have generalized by replacing the element $y$ by an ideal $I$:
this does not affect the proof.
\qed
\end{proof}

The following is similar to Lemma~2.5 in~\cite{shlapentokh2000survey}.

\begin{lemma}
\label{congruencebound}
There exists a constant $c'>0$ 
depending only on $F$ and $K$ such that the following holds:
Let $I$ be a nonzero ideal of $\OO_F$.
Suppose $\mu \in \OO_K$ and $w \in \OO_F$.
Write $\mu = \sum_{i=0}^{s-1} a_i \alpha^i$ with $a_i \in F$.
Suppose $N_{K/\Q}(Da_i) < c' N_{K/\Q}(I)$ for all $i$,
and $\mu \equiv w \pmod{I \OO_K}$.
Then $\mu \in \OO_F$.
\end{lemma}

\begin{proof}
Choose ideals $J_1,\dots,J_h \subseteq \OO_F$ 
representing the elements of the class group of $F$,
and choose $c'>0$ such that $c' N_{K/\Q}(J_j) < 1$ for all $j$.
Choose $j$ such that $J_j I^{-1}$ is principal, 
generated by $z \in F^\times$, say.
Since $\mu \equiv w \pmod{I \OO_K}$,
we have
	$$z(\mu -w) = z(a_0-w) + (za_1) \alpha + \dots 
		+ (z a_{s-1}) \alpha^{s-1} \in \OO_K.$$
By Lemma~4.1 of~\cite{shlapentokh2000survey} (an elementary lemma
about discriminants), $D z a_i \in \OO_F$ 
for $i=1,2,\dots,s-1$.
On the other hand,
	$$|N_{K/\Q}(D z a_i)|  = |N_{K/\Q}(D a_i) N_{K/Q}(z)|
	<  c' N_{K/\Q}(I) \frac{ N_{K/Q}(J_j)}{N_{K/Q}(I)} < 1,$$
by definition of $c'$, so $Dza_i=0$.
Thus $a_i=0$ for $i=1,2,\dots,s-1$.
Hence $\mu \in \OO_F$.
\qed
\end{proof}

\subsection{Denominators of $x$-coordinates of points on an elliptic curve}

We assume that an elliptic curve $E$ as in Theorem~\ref{main} exists.
Thus $E$ is defined over $F$,
and $\rk\, E(F)=\rk\, E(K)=1$.
Hence $E$ has a Weierstrass model of the form $y^2=x^3+ax+b$ 
and we may assume $a,b \in \OO_F$.
Let $O$ denote the point at infinity on $E$,
which is the identity of $E(F)$.

For each nonarchimedean place $\pp$ of $K$, 
let $K_\pp$ denote the completion of $K$ at $\pp$.
and let $\F_\pp$ denote the residue field.
Reducing coefficients modulo $\pp$ yields a possibly singular curve
	$$E_\pp := 
	\Proj \frac{\F_\pp[X,Y,Z]}
		{(Y^2Z - X^3 - \bar{a} XZ^2 - \bar{b} Z^3)}$$
over $\F_\pp$.
Let $E_\pp^\smooth$ denote the smooth part of $E_\pp$.
Let $E_0(K_\pp)$
be the set of points in $E(K_\pp)$ whose reduction mod $\pp$
lies in $E_\pp^\smooth(\F_\pp)$.

\begin{lemma}
\label{E0}
\nichts
\begin{enumerate}
\item $E_0(K_\pp)$ is a subgroup of $E(K_\pp)$.
\item $E_\pp^\smooth(\F_\pp)$ is an abelian group under the usual
chord-tangent law.
\item Reduction modulo $\pp$ gives a surjective group homomorphism 
$\red_\pp: E_0(K_\pp) \rightarrow E_\pp^\smooth(\F_\pp)$.
\item Both $E_0(K_\pp)$ and $E_1(K_\pp):=\ker(\red_\pp)$ are of
finite index in $E(K_\pp)$.
\end{enumerate}
\end{lemma}

\begin{proof}
For the first three statements,
see Proposition~VII.2.1 in~\cite{silvermanAEC}.
We have not assumed that our Weierstrass model is minimal at $\pp$,
so our definition of $E_0$ is different from the standard one
in~\cite{silvermanAEC},
but this does not matter in the proofs.
To prove statement~4,
observe that $E_0(K_\pp)$ and $E_1(K_\pp)$ are open subgroups of 
the compact group $E(K_\pp)$ in the $\pp$-adic topology.
\qed
\end{proof}

{}From now on, $r \in \Z_{\ge 1}$ 
is assumed to be a multiple of $\# E(K)_\tors$,
of the index $(E(K):E(F))$, and of the index $(E(K_\pp):E_0(K_\pp))$ 
for each bad nonarchimedean place $\pp$.
Then $r E(K)$ is a subgroup of $E(F)$ that is free of rank~1,
and $r E(K)$ is contained in $E_0(K_\pp)$ for every $\pp$.

We will need a diophantine approximation result.
First we define the norm $\|\;\|_v: K \rightarrow \R_{\ge 0}$
for each place $v$ of $K$;
it will be characterized by its values on $a \in \OO_K$.
If $v$ is nonarchimedean and $a \in \OO_K - \{0\}$, 
then $\|a\|_v := q^{-v(a)}$
where $q$ is the size of the residue field,
and the discrete valuation $v$ is normalized to take values in $\Z$.
If $v$ is real, then $\|a\|_v$ is the standard absolute value
of the image of $a$ under $K \rightarrow \R$.
If $v$ is complex, then $\|a\|_v$ is the square
of the standard absolute value of the image of $a$ under $K \rightarrow \C$.
Define the naive logarithmic height of $a \in K$ by
	$$h(a):= \sum_{\text{places $v$ of $K$}} \log \max\{\|a\|_v,1\}.$$
If one sums over only the nonarchimedean places $v$,
one obtains $\log N_{K/\Q} \; \den(a)$.

\begin{proposition}
\label{diophantineapproximation}
Let $X$ be a smooth, projective, geometrically integral curve
over $K$ of genus $\ge 1$.
Fix a place $v$ of $K$.
Let $\phi$ be a nonconstant rational function on $X$.
Let $P_1,P_2,\dots$ be a sequence of distinct points in $X(K)$.
For sufficiently large $m$, $P_m$ is not a pole of $\phi$,
so $z_m:=\phi(P_m)$ belongs to $K$.
Then 
	$$\lim_{m \rightarrow \infty} \frac{\log \|z_m\|_v}{h(z_m)} = 0.$$
\end{proposition}

\begin{proof}
See Section~7.4 of~\cite{serremordellweil}.
\qed
\end{proof}

\begin{lemma}
\label{distinctdenominators}
The following holds if $r$ is sufficiently large:
If $P \in r E(K) - \{O\}$ and $m \in \Z - \{-1,0,1\}$,
then 
	$$\log N_{K/\Q} \; \den(x(mP)) \ge \frac{9}{10} m^2 \log N_{K/\Q} \; \den(x(P)) > 0;$$
in particular
$\den(x(mP)) \not= \den(x(P))$
and $\den(x(P)) \not= (1)$.
\end{lemma}

\begin{proof}
Let $P_1$ be a generator of $r E(K)$.
The theory of the canonical height in
Chapter~8, Section~9 of~\cite{silvermanAEC}
implies that there is a real number $\hat{h}(P_1)>0$
(namely, the canonical height of $P_1$, suitably normalized)
such that $h(x(mP_1)) = m^2 \hat{h}(P_1) + O(1)$,
where the implied constant is independent of $m \in \Z$.
Proposition~\ref{diophantineapproximation} applied to each archimedean $v$,
with $X=E$ and $\phi=x$, shows that if we forget to
include the (finitely many) archimedean places
in the sum defining $h$, we obtain
	$$\log N_{K/\Q} \; \den(x(mP_1)) 
	= (1-o(1)) h(x(mP_1)) = (1-o(1)) m^2 \hat{h}(P_1)$$
as $|m| \rightarrow \infty$.
The results follow for large $r$.
\qed
\end{proof}

Of course, there is nothing special about $9/10$;
any real number in the interval $(1/4,1)$ would have done just as well.

\subsection{Divisibility of denominators}

{}From now on, we suppose that $r$ is large enough that
Lemma~\ref{distinctdenominators} holds.

\begin{lemma}
\label{divisibility1}
Let $P,P' \in r E(K) - \{O\}$.
Then $\den(x(P))  \mid  \den(x(P'))$ 
if and only if $P'$ is an integral multiple of $P$.
\end{lemma}

\begin{proof}
We first show that for any ideal $I \subseteq \OO_K$,
the set 
	$$G_I:=\{\, Q \in r E(K) : I  \mid  \den(x(Q)) \,\}$$
is a subgroup of $r E(K)$.
(By convention, we consider $O$ to be an element of $G_I$.)
Since an intersection of subgroups is a subgroup,
it suffices to prove this when $I=\pp^n$
for some prime $\pp$ and some $n \in \Z_{\ge 1}$.
Let $\OO_\pp$ be the completion of $\OO_K$ at $\pp$.
Let $\FF \in \OO_K[[z_1,z_2]]$ 
denote the formal group of $E$ with respect to
the parameter $z:=-x/y$, as in Chapter~4 of~\cite{silvermanAEC}.
Then there is an isomorphism $\FF(\pp \OO_\pp) \isom E_1(K_\pp)$,
given by $z \mapsto (x(z),y(z))$
where $x(z)=z^{-2}+\dots$ and $y(z)=-z^{-3}+\dots$
are Laurent series with coefficients in $\OO_K$.
It follows that $G_{\pp^n}$ is the set of points
in $r E(K)$ lying in the image of $\FF(\pp^{\lceil n/2 \rceil} \OO_\pp)$.
In particular $G_{\pp^n}$ is a subgroup of $r E(K)$.

The ``if'' part of the lemma follows from the preceding paragraph.
Now we prove the ``only if'' part.
Let $G=G_{\den(x(P))}$.
Then $G$ is a subgroup of $r E(K) \isom \Z$,
so $G$ is free of rank~1.
Let $Q$ be a generator of $G$.
By definition of $G$, we have $P \in G$,
so $P$ is a multiple of $Q$.
By the ``if'' part already proved,
\mbox{$\den(x(Q))  \mid  \den(x(P))$}.
On the other hand, $Q \in G$, so \mbox{$\den(x(P)) \mid \den(x(Q))$} 
by definition of $G$.
Thus $\den(x(Q))=\den(x(P))$.
By Lemma~\ref{distinctdenominators}, $Q=\pm P$.
If \mbox{$\den(x(P))  \mid  \den(x(P'))$}, then $P' \in G = \Z Q = \Z P$.
\qed
\end{proof}

\begin{lemma}
\label{anydenominatoroccurs}
If $I \subseteq \OO_K$ is a nonzero ideal,
then there exists $P \in r E(K) - \{O\}$
such that $I \mid \den(x(P))$.
\end{lemma}

\begin{proof}
We use the notation of the previous proof.
It suffices to show that $G_{\pp^n}$ is nontrivial.
This holds since the image of $\FF(\pp^{\lceil n/2 \rceil} \OO_\pp)$
under $\FF(\pp \OO_\pp) \isom E_1(K_\pp)$
is an open subgroup of $E(K_\pp)$, hence of finite index.
\qed
\end{proof}

\begin{lemma}
\label{divisibility2}
Suppose $P \in r E(K) - \{O\}$ and $m \in \Z - \{0\}$.
Let $t=x(P)$ and $t'=x(mP)$.
Then $\den(t)  \mid  \num((t/t' - m^2)^2)$.
\end{lemma}

\begin{proof}
Suppose that $\pp$ is a prime dividing $\den(t)$.
Let $v_\pp: K_\pp \rightarrow \Z \cup \{\infty\}$
denote the discrete valuation associated to~$\pp$.
Then $n:=v_\pp(z(P))$ is positive.
Since $x=z^{-2} + \dots$ is a Laurent series with coefficients in $\OO_K$,
we have $x(P) \in z(P)^{-2} ( 1 + \pp^n \OO_\pp)$.
Using the formal group,
we see that $z(mP) \in m z(P) + \pp^{2n} \OO_\pp$;
in particular $v_\pp(z(mP)) \ge n$,
so $x(mP) \in z(mP)^{-2} ( 1 + \pp^n \OO_\pp)$.
Thus 
	$$\frac{t}{t'} = \frac{x(P)}{x(mP)} \in 
	\left(\frac{z(mP)}{z(P)}\right)^2 \left( 1 + \pp^n \OO_\pp \right).$$
But $\frac{z(mP)}{z(P)} \in m + \pp^n \OO_\pp$,
so $t/t' \in m^2 + \pp^n \OO_\pp$,
so $\pp^n  \mid  \num(t/t'-m^2)$.
On the other hand, $\pp^{2n}$ is the exact power of $\pp$
dividing $\den(t)$.
Applying this argument to every $\pp$
proves \mbox{$\den(t)  \mid  \num((t/t' - m^2)^2)$}.
\qed
\end{proof}

\subsection{Diophantine definition of $\OO_F$ over $\OO_K$}

\begin{lemma}
\label{almostsquares}
With hypotheses as in Theorem~\ref{main},
there exists a subset $S \subseteq \OO_K$ 
such that $S$ is diophantine over $\OO_K$
and $\{\,m^2 : m \in \Z_{\ge 1} \,\} \subseteq S \subseteq \OO_F$.
\end{lemma}

\begin{proof}
Let $c$ and $c'$ be the constants of Lemmas \ref{divisibilitybound} 
and~\ref{congruencebound},
respectively.
By Lemma~\ref{distinctdenominators},
if $\ell \in \Z_{\ge 1}$ is sufficiently large,
then 
	$$c' N_{K/\Q} \den(x(\ell P_0))^{1/2} > N_{K/\Q} \den(x(P_0)^c)$$
for all $P_0 \in r E(K) - \{O\}$.
Fix such an $\ell$.

Let $S$ be the set of $\mu \in \OO_K$
such that there exist $P_0,P',P' \in r E(K) - \{O\}$
and $t_0,t,t' \in F$ such that 
\begin{enumerate}
\item $P=\ell P_0$
\item $t_0=x(P_0)$, $t=x(P)$, $t'=x(P')$
\item $(\mu+1)(\mu+2)\dots(\mu+n) \mid \den(t_0)$
\item $\den(t) \mid \den(t')$
\item $\den(t) \mid \num((t/t' - \mu)^2)$
\end{enumerate}
It follows from Lemma~\ref{ideals} that $S$ is diophantine over $\OO_K$.

Suppose $m \in \Z_{\ge 1}$.
We wish to show that $\mu:=m^2$ belongs to $S$.
By Lemma~\ref{anydenominatoroccurs},
there exists $P_0 \in r E(K) - \{O\}$
such that $(\mu+1)(\mu+2)\dots(\mu+n) \mid \den(x(P_0))$.
Let $P=\ell P_0$ and $P'=mP$.
Let $t_0=x(P_0)$, $t=x(P)$, and $t'=x(P')$.
Then conditions (1), (2), and~(3) in the definition of $S$
are satisfied,
and (4) and~(5) follow from Lemmas \ref{divisibility1}
and~\ref{divisibility2}, respectively.
Hence $m^2 \in S$.

Now suppose that $\mu \in S$.
We wish to show that $\mu \in \OO_F$.
Fix $P_0$, $P$, $P'$, $t_0$, $t$, $t'$
satisfying~(1) through~(5).
By~(4) and Lemma~\ref{divisibility1}, $P'=mP$ for some nonzero $m \in \Z$.
By Lemma~\ref{divisibility2}, 
\mbox{$\den(t) \mid \num((t/t'-m^2)^2)$}.
On the other hand, (5) says that \mbox{$\den(t) \mid \num((t/t'-\mu)^2)$}.
Therefore $\den(t)^{1/2} \mid \num(\mu - m^2) = (\mu-m^2)$.
(Note that each prime of $\OO_F$ or of $\OO_K$ that appears in
$\den(t)$ must occur to an even power, since $t$ is the $x$-coordinate
of a point on $y^2=x^3+ax+b$.  Hence $\den(t)^{1/2}$ is a well-defined
ideal.)
Write $\mu=\sum_{i=0}^{s-1} a_i \alpha^i$ with $a_i \in F$.
By~(3) and Lemma~\ref{divisibilitybound}, 
$N_{K/\Q}(Da_i) \le N_{K/\Q}(\den(t_0))^c$.
By definition of $\ell$, we have
$N_{K/\Q}(\den(t_0))^c < c' N_{K/\Q} \den(t)^{1/2}$.
Combining these shows
that the hypotheses of Lemma~\ref{congruencebound}
hold for $w=m^2$ and $I=\den(t)^{1/2}$ (as an ideal in $\OO_F$).
Thus $\mu \in \OO_F$.
\qed
\end{proof}

\noindent
{\em Proof of Theorem~\ref{main}.}
Let $S$ be the set given by Lemma~\ref{almostsquares}.
Then $S_1:=\{\, s-s' : s,s' \in S \,\}$
contains all odd integers at least $3$,
because of the identity $(m+1)^2-m^2=2m+1$.
Next, $S_2:=S_1 \cup \{\, 4-s: s \in S_1 \,\}$
contains all odd integers,
and $S_3:=S_2 \cup \{\, s+1: s \in S_2 \,\}$
contains $\Z$.
Let $\beta_1,\dots,\beta_b$ be a $\Z$-basis for $\OO_F$.
Then 
$S_4:= \{\, a_1 \beta_1 + \dots + a_b \beta_b : a_1,\dots,a_b \in S_3 \,\}$
contains $\OO_F$.

But $S \subseteq \OO_F$, so $S_i \subseteq \OO_F$ for $i=1,2,3,4$.
In particular, $S_4=\OO_F$.
Also, $S$ is diophantine over $\OO_K$,
so each $S_i$ is diophantine over $\OO_K$.
In particular, $\OO_F=S_4$ is diophantine over $\OO_K$.
\qed

\subsection{Questions}

\begin{enumerate}
\item Is it true that for every number field $K$,
there exists an elliptic curve $E$ over $\Q$
such that $\rk\, E(\Q) = \rk\, E(K) = 1$?
The author would conjecture so.
If so, then Hilbert's Tenth Problem over $\OO_K$
is undecidable for every number field $K$.

\item Can one weaken the hypotheses of Theorem~\ref{main}
and give a diophantine definition of $\OO_F$ over $\OO_K$
using any elliptic curve $E$ over $K$ with $\rk\, E(K)=1$,
not necessarily defined over $F$?
Such elliptic curves may be easier to find.
But our proof of Theorem~\ref{main}
seems to require the fact that $E$ is defined over $F$ and has $\rk\, E(F)=1$,
since Lemma~\ref{congruencebound} fails 
if the ideal $I$ of $\OO_F$ is instead assumed to be an ideal of $\OO_K$.

\item Can one prove an analogue of Theorem~\ref{main}
in which the elliptic curve is replaced by an abelian variety?
\end{enumerate}


\providecommand{\bysame}{\leavevmode\hbox to3em{\hrulefill}\thinspace}


\begin{thebibliography}{DLPVG00}

\bibitem[CZ00]{cornelissen-zahidi2000}
Gunther Cornelissen and Karim Zahidi, \emph{Topology of {D}iophantine sets:
  remarks on {M}azur's conjectures}, Hilbert's tenth problem: relations with
  arithmetic and algebraic geometry (Ghent, 1999), Amer. Math. Soc.,
  Providence, RI, 2000, pp.~253--260.

\bibitem[Dav53]{davis1953}
Martin Davis, \emph{Arithmetical problems and recursively enumerable
  predicates}, J. Symbolic Logic \textbf{18} (1953), 33--41.

\bibitem[Den80]{denef1980}
J.~Denef, \emph{Diophantine sets over algebraic integer rings. {I}{I}}, Trans.
  Amer. Math. Soc. \textbf{257} (1980), no.~1, 227--236.

\bibitem[DL78]{denef-lipshitz1978}
J.~Denef and L.~Lipshitz, \emph{Diophantine sets over some rings of algebraic
  integers}, J. London Math. Soc. (2) \textbf{18} (1978), no.~3, 385--391.

\bibitem[DLPVG00]{hilbertstenthproblem}
Jan Denef, Leonard Lipshitz, Thanases Pheidas, and Jan Van~Geel (eds.),
  \emph{Hilbert's tenth problem: relations with arithmetic and algebraic
  geometry}, American Mathematical Society, Providence, RI, 2000, Papers from
  the workshop held at Ghent University, Ghent, November 2--5, 1999.

\bibitem[DPR61]{davis-putnam-robinson1961}
Martin Davis, Hilary Putnam, and Julia Robinson, \emph{The decision problem for
  exponential diophantine equations}, Ann. of Math. (2) \textbf{74} (1961),
  425--436.

\bibitem[Eis]{eisentraeger2003}
Kirsten Eisentr\"ager, Ph.~D.\ thesis, University of California, Berkeley, in
  preparation.

\bibitem[KR92]{kim-roush1992}
K.~H. Kim and F.~W. Roush, \emph{Diophantine undecidability of {${\bf {C}}(t\sb
  1,t\sb 2)$}}, J. Algebra \textbf{150} (1992), no.~1, 35--44.

\bibitem[Mat70]{matijasevic}
Ju.~V. Matijasevi{\v{c}}, \emph{The {D}iophantineness of enumerable sets},
  Dokl. Akad. Nauk SSSR \textbf{191} (1970), 279--282.

\bibitem[Maz94]{mazur1994}
B.~Mazur, \emph{Questions of decidability and undecidability in number theory},
  J. Symbolic Logic \textbf{59} (1994), no.~2, 353--371.

\bibitem[MB]{moretbailly2002}
Laurent Moret-Bailly, paper in preparation, extending results presented in a
  lecture 18 June 2001 at a conference in honor of Michel Raynaud in Orsay,
  France.

\bibitem[Phe88]{pheidas1988}
Thanases Pheidas, \emph{Hilbert's tenth problem for a class of rings of
  algebraic integers}, Proc. Amer. Math. Soc. \textbf{104} (1988), no.~2,
  611--620.

\bibitem[Phe91]{pheidas1991}
Thanases Pheidas, \emph{Hilbert's tenth problem for fields of rational
  functions over finite fields}, Invent. Math. \textbf{103} (1991), no.~1,
  1--8.

\bibitem[Phe00]{pheidas2000}
Thanases Pheidas, \emph{An effort to prove that the existential theory of
  {${\bf {Q}}$} is undecidable}, Hilbert's tenth problem: relations with
  arithmetic and algebraic geometry (Ghent, 1999), Amer. Math. Soc.,
  Providence, RI, 2000, pp.~237--252.

\bibitem[PZ00]{pheidas-zahidi2000}
Thanases Pheidas and Karim Zahidi, \emph{Undecidability of existential theories
  of rings and fields: a survey}, Hilbert's tenth problem: relations with
  arithmetic and algebraic geometry (Ghent, 1999), Amer. Math. Soc.,
  Providence, RI, 2000, pp.~49--105.

\bibitem[Ser97]{serremordellweil}
Jean-Pierre Serre, \emph{Lectures on the {M}ordell-{W}eil theorem}, third ed.,
  Friedr. Vieweg \& Sohn, Braunschweig, 1997, Translated from the French and
  edited by Martin Brown from notes by Michel Waldschmidt, With a foreword by
  Brown and Serre.

\bibitem[Shl89]{shlapentokh1989}
Alexandra Shlapentokh, \emph{Extension of {H}ilbert's tenth problem to some
  algebraic number fields}, Comm. Pure Appl. Math. \textbf{42} (1989), no.~7,
  939--962.

\bibitem[Shl92]{shlapentokh1992functionfield}
Alexandra Shlapentokh, \emph{Hilbert's tenth problem for rings of algebraic
  functions in one variable over fields of constants of positive
  characteristic}, Trans. Amer. Math. Soc. \textbf{333} (1992), no.~1,
  275--298.

\bibitem[Shl00a]{shlapentokh2000functionfield}
Alexandra Shlapentokh, \emph{Hilbert's tenth problem for algebraic function
  fields over infinite fields of constants of positive characteristic}, Pacific
  J. Math. \textbf{193} (2000), no.~2, 463--500.

\bibitem[Shl00b]{shlapentokh2000survey}
Alexandra Shlapentokh, \emph{Hilbert's tenth problem over number fields, a
  survey}, Hilbert's tenth problem: relations with arithmetic and algebraic
  geometry (Ghent, 1999), Amer. Math. Soc., Providence, RI, 2000, pp.~107--137.

\bibitem[Sil92]{silvermanAEC}
Joseph~H. Silverman, \emph{The arithmetic of elliptic curves}, Springer-Verlag,
  New York, 1992, Corrected reprint of the 1986 original.

\bibitem[Vid94]{videla1994}
Carlos~R. Videla, \emph{Hilbert's tenth problem for rational function fields in
  characteristic $2$}, Proc. Amer. Math. Soc. \textbf{120} (1994), no.~1,
  249--253.

\end{thebibliography}
\end{document}